\documentclass[11pt]{article}
\usepackage{}

\usepackage[paperwidth=9.5in,top=1.3in, bottom=1.4in]{geometry}
\usepackage{amsfonts}
\usepackage{amsthm}
\usepackage{amssymb}
\usepackage{latexsym}
\usepackage{amsmath,amsfonts}
\usepackage{mathrsfs}
\usepackage{cases}
\usepackage{latexsym,bm}
\usepackage{indentfirst}
\usepackage{color}
\usepackage{ifpdf}
\usepackage{graphicx}
\usepackage{psfrag}
\usepackage{dsfont}
\usepackage[
pdfauthor={CSY},
pdftitle={Strong generalized Halin graphs},
pdfstartview=XYZ,
bookmarks=true,
colorlinks=true,
linkcolor=blue,
urlcolor=blue,
citecolor=blue,
bookmarks=true,
linktocpage=true,
hyperindex=true
]{hyperref}

\newcommand{\xiaowuhao}{\fontsize{9pt}{\baselineskip}\selectfont}

\newtheorem{THM}{\textbf{Theorem}}[section]
\newtheorem{DEF}{\textbf{Definition}}
\newtheorem{LEM}{\textbf{Lemma}}[section]
\newtheorem{CON}{\textbf{Conjecture}}
\newtheorem{COR}{\textbf{Corollary}}[section]

\newcommand{\pf}{\textbf{Proof}.\quad}

\newtheorem{CLA}{\textbf{Claim}}[section]
\newcommand{\phibar}{\overline{\varphi}}
\newcommand{\phiv}{\varphi}
\begin{document}

\title{Chromatic index determined by fractional chromatic index}

\author{Guantao Chen$^{a}$, Yuping Gao$^{b,a}$, Ringi Kim$^{c}$,  Luke Postle$^c$,  Songling Shan$^{d}$\\
{\xiaowuhao  $^{a}$ Department of Mathematics and Statistics, Georgia State University, Atlanta, GA\,30303, USA}\\
{\xiaowuhao $^{b}$ School of Mathematics and Statistics, Lanzhou University, Lanzhou, Gansu 730000, China}\\
{\xiaowuhao $^c$  University of Waterloo,  Waterloo, ON, N2L 3G1, Canada}\\
{\xiaowuhao $^{d}$ Department of Mathematics, Vanderbilt University, Nashville, TN\, 37240, USA}}

\date{}

\maketitle

\begin{abstract}Given a graph $G$ possibly with multiple edges but no loops, denote by $\Delta$ the maximum degree, $\mu$ the multiplicity, $\chi'$ the chromatic index and $\chi_f'$ the fractional chromatic index of $G$, respectively.   It is known that $\Delta\le \chi_f' \le \chi' \le \Delta + \mu$, where the upper bound is a classic result of Vizing.  While deciding the exact value of $\chi'$ is a classic NP-complete problem, the computing of $\chi_f'$ is in polynomial time. In fact, it is shown that  if $\chi_f' > \Delta$ then $\chi_f'= \max  \frac{|E(H)|}{\lfloor |V(H)|/2\rfloor}$, where the maximality is taken over all induced subgraphs $H$ of $G$.  Gupta\,(1967), Goldberg\,(1973), Andersen\,(1977), and Seymour\,(1979) conjectured that  $\chi'=\lceil\chi_f'\rceil$ if $\chi'\ge \Delta+2$, which is commonly referred as Goldberg's conjecture. It has been shown that Goldberg's conjecture is equivalent to the following conjecture of Jakobsen:  For any positive integer $m$ with $m\ge 3$,   every graph $G$ with $\chi'>\frac{m}{m-1}\Delta+\frac{m-3}{m-1}$ satisfies $\chi'=\lceil\chi_f'\rceil$. Jakobsen's conjecture has been verified for $m$ up to 15 by various researchers in the last four decades. We use an extended form of a Tashkinov tree to show that it is true  for $m\le 23$. With the same technique, we show that if $\chi' \geq\Delta+\sqrt[3]{\Delta/2}$ then $\chi'=\lceil\chi_f'\rceil$. The previous best known result is for graphs with $\chi'> \Delta +\sqrt{\Delta/2}$ obtained by Scheide, and by Chen, Yu and Zang, independently.
  Moreover, we show that Goldberg's conjecture holds for graphs $G$ with $\Delta\leq 23$ or $|V(G)|\leq 23$.
\end{abstract}



\section{Introduction}

Graphs considered in this paper may contain multiple edges but no loops. Let
$G$ be a graph and $\Delta:=\Delta(G)$ be the maximum degree of $G$.
A (proper) {\it  $k$-edge-coloring} $\varphi$ of $G$ is a mapping  $\varphi$ from $E(G)$ to $\{1, 2, \cdots, k\}$ (whose elements are called colors)
such that no two adjacent edges receive the same color.
The {\it chromatic index}  $\chi' :=\chi'(G)$
is  the least integer $k$ such that $G$ has a $k$-edge-coloring.
In graph edge-coloring, the central question is to determine the chromatic index $\chi'$ for graphs. We refer the book~\cite{StiebSTF-Book} of Stiebitz, Scheide, Toft and Favrholdt and the elegant survey~\cite{McDonaldSurvey15} of McDonald for literature on the recent progress of graph edge-colorings.
Clearly, $\chi' \ge \Delta$.  Conversely, Vizing showed that $\chi' \le \Delta + \mu$, where $\mu := \mu(G)$ is the multiplicity of $G$.  However, determining the exact value of $\chi'$ is a very difficult problem.  Holyer~\cite{Holyer81} showed that the problem is NP-hard even restricted to simple cubic  graphs.  To estimate $\chi'$, the notion of fractional chromatic index is introduced.

A {\it fractional edge coloring} of $G$ is a non-negative weighting $w(.)$ of the set $\mathcal{M}(G)$ of
matchings in $G$ such that, for every edge $e\in E(G)$, $\sum_{M\in \mathcal{M}: e\in M} w(M) =1$.  Clearly, such a weighting $w(.)$ exists.  The {\it fractional chromatic index} $\chi_f' := \chi'_f(G)$ is the  minimum total weight $\sum_{M\in \mathcal{M}} w(M)$ over all fractional edge colorings of $G$. By definitions, we have $\chi' \ge \chi'_f\ge \Delta $.
It follows from Edmonds' characterization of the matching polytope~\cite{Edmonds65} that $\chi_f'$ can be computed in polynomial time and
\[
\chi'_f = \max\left\{ \frac{|E(H)|}{\lfloor |V(H)|/2\rfloor} \, : \,   \mbox{  $H\subseteq G$ with $|V(H)|\ge 3$ } \right\} \, \mbox{if $\chi_f' > \Delta$}.
\]
It is not difficult to show that the above maximality can be restricted to induced subgraphs $H$ with odd number of vertices.  So, in the case of $\chi'_f > \Delta$, we have
$$\lceil \chi_f'\rceil=\max\left\{\left\lceil \frac{2|E(H)|}{|V(H)|-1}\right\rceil:  \mbox{induced\ subgraphs\ $H\subseteq G$\ with\ $|V(H)| \ge 3$ \ and\ odd} \,\right\}.$$

A graph $G$ is called {\it elementary} if $\chi' = \lceil \chi_f' \rceil$.   Gupta\,(1967)~\cite{Gupta67}, Goldberg\,(1973)~\cite{Goldberg}, Andersen\,(1977)~\cite{Andersen}, and Seymour\,(1979)~\cite{Seymour} independently made the following conjecture, which is commonly referred as {\it Goldberg's conjecture}.

\begin{CON}\label{con:GAS}
For any graph $G$,  if $\chi' \ge \Delta+2$ then  $G$ is elementary.
\end{CON}

An immediate consequence of Conjecture~\ref{con:GAS} is that $\chi'$ can be computed in polynomial time for graphs with $\chi' \ge \Delta +2$. So the NP-complete problem of computing the chromatic indices lies in determining whether $\chi' = \Delta$,  $\Delta +1$, or $\ge \Delta +2$, which strengthens Vizing's classic result  $\chi'\leq \Delta+\mu$ tremendously when $\mu$ is big.

 Following $\chi'\le \frac{3\Delta}2$ of the classic result of Shannon~\cite{Shannon49}, we can assume that, for every $\Delta$,  there exists the least positive number $\zeta$ such that if $\chi'> \Delta + \zeta$ then $G$ is elementary.  Conjecture~\ref{con:GAS} indicates that $\zeta \le 1$. Asymptotically, Kahn~\cite{Kahn96} showed $\zeta= o(\Delta)$. Scheide~\cite{Scheide-2010}, and Chen, Yu,  and Zang~\cite{CYZ-2011} independently proved that
$\zeta \le \sqrt{\Delta/2}$.  In this paper, we show that $\zeta \le \sqrt[3]{\Delta/2}-1$ as stated below.

\begin{THM}\label{THM:cubic}
For any graph $G$,   if   $\chi'\geq\Delta+\sqrt[3]{\Delta/2}$,
then $G$ is  elementary.
\end{THM}

Jakobsen~\cite{Jakobsen73} conjectured that $\zeta \le 1 + \frac{\Delta -2}{m-1}$ for every positive integer $m(\geq3)$, which gives a reformulation of Conjecture~\ref{con:GAS} as stated below.
\begin{CON}\label{con:Jm}
Let $m$  be an integer with $m\ge 3$ and $G$ be a graph.  If
$\chi'>\frac{m}{m-1}\Delta+\frac{m-3}{m-1}$, then $G$ is elementary.
\end{CON}
Since $\frac{m}{m-1}\Delta+\frac{m-3}{m-1}$ decreases as $m$ increases,  it is sufficient to prove Jakobsen's conjecture for all odd integers $m$ (in fact, for any infinite sequence of positive integers), which has been confirmed slowly for $m\le 15$ by a series of papers over the last 40 years:
\begin{itemize}
\item {\bf $m=5$:}   Three independent proofs given by Andersen~\cite{Andersen} (1977), Goldberg~\cite{Goldberg} (1973),  and S{\o}rensen\,(unpublished, page 158 in \cite{StiebSTF-Book}), respectively.

\item {\bf $m=7$:}  Two independent proofs given by Andersen~\cite{Andersen} (1977)
and S{\o}rensen (unpublished, page 158 in \cite{StiebSTF-Book}), respectively.

\item {\bf $m=9$:} By  Goldberg~\cite{Goldberg-1984} (1984).

\item {\bf $m=11$:}   Two independent proofs given by Nishizeki and Kashiwagi~\cite{Nishizeki-Kashiwagi-1990} (1990) and by Tashkinov~\cite{Tashkinov-2000} (2000), respectively.

\item {\bf $m=13$:}   By Favrholdt, Stiebitz and Toft~\cite{FavrST06} (2006).

\item {\bf $m=15$:} By   Scheide~\cite{Scheide-2010} (2010).
\end{itemize}

In this paper, we show that Jakobsen's conjecture is true up to $m=23$.
\begin{THM}\label{THM:Jm19}
If  $G$ is a graph with $\chi'>\frac{23}{22}\Delta+\frac{20}{22}$,
then $G$ is elementary.
\end{THM}

\begin{COR}\label{COR:Delta23}
If  $G$ is a graph with $\Delta\le 23$ or $|V(G)|\le 23$, then $\chi'\le \max\{\Delta+1,\lceil\chi'_f\rceil\}$.
\end{COR}

Note that in Corollary~\ref{COR:Delta23}, $|V(G)| \le 23$ does not imply $\Delta \le 23$, as $G$ may have multiple edges.
The remainder of this paper is organized as follows. In Section 2, we introduce some definitions and notation
for edge-colorings, Tashkinov trees, and several known results which are useful for the proofs of Theorems~\ref{THM:cubic} and~\ref{THM:Jm19};
in Section 3, we give an extension of Tashkinov trees and prove several properties of the extended Tashkinov trees; and
in Section 4, we prove Theorem~\ref{THM:cubic}, Theorem~\ref{THM:Jm19} and Corollary \ref{COR:Delta23} based on the results in Section 3.

\section{Preliminaries}

\subsection{Basic definitions and notation}

Let $G$ be a graph with vertex set $V$ and edge set $E$.  Denote by $|G|$ and $||G||$ the number of vertices and the number of edges of  $G$, respectively.  For any two sets $X, Y\subseteq V$,  denote by  $E(X, Y)$ the set of edges with one end in $X$ and the other one in $Y$ and
denote by $\partial(X) := E(X, V-X)$  the  boundary edge set of $X$, that is, the set of edges with exactly one end in $X$.  Moreover,  let $E(x, y) := E(\{x\}, \{y\})$  and  $E(x) := \partial(\{x\})$.   Denote by $G[X]$ the subgraph induced by $X$ and $G-X$ the subgraph induced by $V(G)-X$.
Moreover,  let $G-x = G-\{x\}$.   For any subgraph $H$ of $G$, we let $G[H] = G[V(H)]$ and $\partial(H) = \partial(V(H))$. Let $V(e)$ be the set of the two ends of an edge $e$.

A path $P$ is usually denoted by an alternating sequence $P=(v_0, e_1, v_1,\cdots, e_p,v_p)$ with $V(P)=\{v_0,\cdots, v_p\}$
and $E(P)=\{e_1,\cdots, e_p\}$ such that $e_i\in E_G(v_{i-1}, v_i)$
for $1\le i\le p$. The path $P$ defined above  is called
 a $(v_0,v_p)$-path.  For any two vertices  $u,v \in V(P)$,
denote by $uPv$ or $vPu$ the unique subpath connecting $u$ and $v$. If $u$ is an end of $P$, then
we obtain a {\it linear order  $\preceq_{(u, P)}$} of the vertices of $P$ in a natural way such that  $x\preceq_{(u,P)}y$ if  $x\in V(uPy)$.

The set of all $k$-edge-colorings of a graph $G$
is denoted by $\mathcal{C}^k(G)$. Let $\varphi\in \mathcal{C}^k(G)$.  For any color $\alpha$, let
$E_{\alpha} =\{e\in E \ : \  \phiv(e) =\alpha\}$. More generally, for each subgraph $H\subseteq G$, let $$E_{\alpha}(H)=\{e\in E(H)\ :\ \phiv(e)=\alpha\}.$$  For any two distinct colors $\alpha$ and $\beta$, denote
by $G_{\varphi}(\alpha,\beta)$
 the subgraph of $G$ induced by $E_{\alpha} \cup E_{\beta}$. The components of $G_{\varphi}(\alpha,\beta)$ are called
{\it $(\alpha,\beta)$-chains}.  Clearly, each $(\alpha,\beta)$-chain is either a path or a cycle of edges  alternately colored with $\alpha$ and $\beta$.
For each $(\alpha, \beta)$-chain $P$, let $\varphi/P$ denote the $k$-edge-coloring obtained from $\varphi$ by exchanging colors $\alpha$ and $\beta$ on $P$, that is, for each $e\in E$,
$$
\varphi/P\ (e) =
\left\{
  \begin{array}{ll}
    \varphi(e), & \hbox{$e\notin E(P)$;} \\
    \beta, & \hbox{$e\in E(P)$ and $\varphi(e)=\alpha$;} \\
    \alpha, & \hbox{$e\in E(P)$ and $\varphi(e)=\beta$.}
  \end{array}
\right.
$$

For any $ v\in V$, let $P_v(\alpha,\beta,\varphi)$ denote
the unique $(\alpha,\beta)$-chain containing $v$.   Notice that, for any two vertices $u, \, v\in V$, either $P_u(\alpha,\beta,\varphi)=P_v(\alpha,\beta,\varphi)$
or $P_u(\alpha,\beta,\varphi)\cap P_v(\alpha,\beta,\varphi)=\emptyset$.
For any $v\in V$, let $\varphi(v) :=\{\varphi(e)\,: e\in E(v)\}$ denote the set of colors presented at $v$  and
$  \phibar(v)$ the set of colors not assigned to any edge incident to $v$, which are called  {\it missing} colors at $v$.  For any vertex set  $X\subseteq V$, let $\varphi (X) = \cup_{x\in X} \varphi(x)$ and $\phibar(X) = \cup_{x\in X} \phibar(x)$ be the set of colors presenting and missing at some vertices of $X$, respectively. For any edge set $F\subseteq E$, let $\varphi (F) = \cup_{e\in F} \varphi(e)$.

\subsection{Elementary sets and closed sets}

Let $G$ be a graph.  An edge $e\in E(G)$ is called {\it critical} if $\chi'(G-e) < \chi'(G)$,  and
the graph $G$ is called {\it critical} if $\chi'(H) < \chi'(G)$ for any proper subgraph $H\subseteq G$. A graph $G$ is called {\it $k$-critical} if it is critical and $\chi'(G) = k+1$.  In the proofs, we will consider a graph $G$ with $\chi'(G) = k+1 \ge \Delta +2$, a critical edge $e\in E(G)$, and a coloring $\phiv\in \mathcal{C}^k(G-e)$. We call them together a {\it $k$-triple} $(G, e, \phiv)$.

\begin{DEF}
Let $G$ be a graph and $e\in E(G)$ such that $\mathcal{C}^k(G-e)\ne \emptyset$ and let $\varphi \in \mathcal{C}^k(G-e)$. Let $X \subseteq V(G)$ contain two ends of $e$.
\begin{itemize}
\item We call $X$  {\it elementary} \emph{(}with respect to $\varphi$\emph{)} if all missing color sets $\phibar(x)$ \emph{(}$x\in X$\emph{)} are mutually  disjoint.
\item We call $X$ {\it closed}  \emph{(}with respect to $\varphi$\emph{)} if $\varphi(\partial(X))\cap \phibar(X)=\emptyset$, i.e., no missing color of  $X$  appears on the edges in $\partial(X)$. If additionally, each color in $\varphi(X)$ appears at most once in $\partial(X)$,  we call $X$  {\it strongly closed} \emph{(}with respect to
$\varphi$\emph{)}.
\end{itemize}
\end{DEF}

Moreover,  we call a subgraph $H\subseteq G$ {\it elementary}, {\it closed}, and {\it strongly closed}
if  $V(H)$ is elementary, closed, and {strongly closed}, respectively. If a vertex set $X\subseteq V(G)$ containing two ends of $e$ is both elementary and strongly closed, then $|X|$ is odd and $k= \frac{2(|E(G[X])|-1)}{|X| -1}$, so $k +1=\left\lceil\frac{2|E(G[X])|}{|X| -1}\right\rceil=\lceil \chi_f' \rceil$. Therefore, if $V(G)$ is elementary then $G$ is elementary, i.e., $\chi'(G) =k+1 =\lceil \chi_f'\rceil$.

\subsection{Tashkinov trees}

\begin{DEF}\label{Def:Tashkinov-tree}
 A {\it Tashkinov tree} of a $k$-triple $(G, e, \phiv)$ is a tree $T$, denoted by $T=(e_1, e_2, \cdots, e_p)$,
  induced by a sequence of edges $e_1=e$, $e_2$, $\dots$, $e_p$ such that  for each $i\ge 2$, $e_i$ is a boundary edge of the tree induced by $\{e_1, e_2, \cdots, e_{i-1}\}$ and $\varphi(e_i)\in \phibar\left(V\left(\bigcup\limits_{j=1}^{i-1}e_j\right)\right)$.
\end{DEF}

  For each  $e_j\in
\{e_1,\cdots, e_p\}$, we denote by $Te_j$  the subtree
$T[\{e_1,\cdots, e_j\}]$ and denote by $e_jT$ the subgraph induced by  $\{e_j,\cdots,e_p\}$. For each edge $e_i$ with $i \ge 2$, the end of $e_i$ in $Te_{i-1}$ is called the {\it in-end} of $e_i$ and the other one is called the {\it out-end} of $e_i$.

Algorithmically, a Tashkinov tree is obtained incrementally from $e$ by adding a boundary edge whose color is missing in the previous tree.  Vizing-fans (stars) (used in the proof of Vizing's classic theorem~\cite{Vizing64})
 and Kierstead-paths (used in ~\cite{Kierstead84}) are special Tashkinov trees.

 \begin{THM}\label{THM:TashOrigi}$[$Tashkinov $\cite{Tashkinov-2000}$ $]$
For any given $k$-triple $(G, e, \phiv)$ with $k\geq \Delta+1$, all Tashkinov trees are elementary.
\end{THM}

For a graph $G$,  a Tashkinov tree is associated with an edge $e\in E(G)$ and a $k$-edge-coloring of $G-e$ with $k \ge \Delta +1$.  We distinguish  the following three different types of maximality.
\begin{DEF}\label{LEM:T-Property}
Let $(G,e, \phiv)$ be a $k$-triple with $k \ge \Delta +1$, and $T$ be a Tashkinov tree of $(G,e,\phiv)$.
\begin{itemize}
\item We call  $T$  {\it $(e,\phiv)$-maximal} if there is no Tashkinov tree $T^*$  of $(G, e, \phiv)$ containing $T$ as a proper subtree, and denote by $\mathcal{T}_{e, \phiv}$  the set of all $(e,\phiv)$-maximal Tashkinov trees.
\item We call $T$  {\it $e$-maximal} if there is no Tashkinov tree $T^*$ of a $k$-triple $(G, e, \phiv^*)$ containing $T$ as a proper subtree, and  denote by $\mathcal{T}_e$  the set of all $e$-maximal Tashkinov trees.

\item We call $T$ {\it maximum} if $|T|$ is maximum over  all Tashkinov trees of $G$, and denote by $\mathcal{T}$  the set of all maximum Tashkinov trees.
\end{itemize}
\end{DEF}

 Let $T$ be a Tashkinov tree of a $k$-triple $(G, e, \phiv)$.
Then, $T$ is $(e, \phiv)$-maximal if and only if  $V(T)$ is closed. Moreover,
 the  vertex sets are the same for all $T\in \mathcal{T}_{e, \phiv}$.
   We call colors in  $\phiv(E(T))$
{\it used} and colors  not in $\phiv(E(T))$  {\it   unused} on $T$,  call an unused missing color in
$\phibar(V(T))$  a {\it free color} of $T$  and denote the set of all free colors of $T$ by $\Gamma^f(T)$. For each color $\alpha$, let
 $E_{\alpha}(\partial(T))$ denote the set of edges with color $\alpha$ in boundary $\partial(T)$.  A color $\alpha$ is called a {\it defective color} of $T$  if $|E_{\alpha}(\partial(T))|\ge  2$. The set of
all defective colors of $T$ is denoted by $\Gamma^d(T)$.  Note that if $T\in \mathcal{T}_{e, \phiv}$, then
$V(T)$ is strongly closed if and only if $T$ does not have any defective colors.

The following corollary  follows immediately from the fact that a maximal Tashkinov tree is elementary and closed.
\begin{COR}\label{COR-(e,phi)Max}
For each $T\in \mathcal{T}_{e, \phiv}$,  the following properties hold.
\item[$(1)$] $|T| \ge 3$ is odd.
 \item [$(2)$] For any two missing colors $\alpha,\beta\in \phibar(V(T))$, we have $P_u(\alpha,\beta,\phiv)=P_v(\alpha,\beta,\phiv)$, where $u$ and $v$ are the two unique vertices in $V(T)$ such that  $\alpha\in \phibar(u)$ and $\beta\in \phibar(v)$, respectively. Furthermore, $V(P_u(\alpha,\beta,\phiv))\subseteq V(T)$.
 \item[$(3)$] For every defective color $\delta\in \Gamma^d (T)$,  $|E_{\delta}(\partial(T))|\ge 3$ and is odd.
 \item [$(4)$] There are at least four free colors. More specifically,  $$|\Gamma^f(T)|\ge |T|(k-\Delta)+2-|\phiv(E(T))|\ge |T|+2-(|T|-2) \ge 4.$$
  \end{COR}

  The following lemma was given
in~\cite{StiebSTF-Book}.
\begin{LEM}\label{LEM:PassAll}
Let $T\in \mathcal{T}_{e}$ be a Tashkinov tree of a $k$-triple $(G, e, \phiv)$ with $k\geq \Delta+1$.  For any free color $\gamma\in \Gamma^f(T)$ and any $\delta \notin \phibar(V(T))$, the $(\gamma, \delta)$-chain $P_u(\gamma, \delta, \phiv)$ contains all edges in $E_{\delta}(\partial(T))$, where $u$ is the unique vertex of $T$ missing color $\gamma$.

\end{LEM}
\proof Otherwise, consider the coloring $\varphi_1 = \phiv/P_u(\gamma, \delta, \phiv)$. Since $\delta$ and $\gamma$ are both unused on $T$ with respect to $\varphi$, $T$ is still a Tashkinov tree and $\delta$ is a missing color with respect to $\varphi_1$. But $E_{\delta}(\partial(T))\neq \emptyset$, which gives a contradiction to  $T$ being an $e$-maximal tree. \qed

Following the notation in Lemma~\ref{LEM:PassAll}, we consider the case of $\delta$ being a defective
color.   Then $P: =P_u(\gamma,\delta,\phiv)$ is a path with $u$ as one end.
Since $u$ is the unique vertex in $T$ missing $\gamma$ by Theorem~\ref{THM:TashOrigi},
the other end of $P$ is not in $T$.
In the linear order $\preceq_{(u,P)}$, the last vertex $v$ with $v\in V(T)\cap V(P)$ is called an {\it exit vertex} of $T$.
Applying Lemma~\ref{LEM:PassAll}, Scheide~\cite{Scheide-2010} obtained the following result.
\begin{LEM}\label{LEM:exit-color}
Let $T\in \mathcal{T}_e$ be a Tashkinov tree of a $k$-triple $(G, e,\phiv)$ with $k \ge \Delta +1$.
If  $v$ is an exit vertex of $T$, then every missing color in $\phibar(v)$ must be used
on $T$.
\end{LEM}

Let $T\in \mathcal{T}_{e, \phiv}$ be a Tashkinov tree of $(G, e, \phiv)$ and $V(e) =\{x, y\}$.
By keeping odd number of vertices in each step of growing a Tashkinov tree from $e$, Scheide~\cite{Scheide-2010} showed that there is another $T^*\in \mathcal{T}_{e, \phiv}$, named  a {\it  balanced Tashkinov tree},  such that
$V(T^*) = V(T)$ constructed incrementally from $e$ by the following steps:
\begin{itemize}
\item {\bf Adding a path:} Pick two missing colors $\alpha$ and $\beta$ with $\alpha \in \phibar(x)$ and $\beta\in \phibar(y)$,  and let $T^*:= \{e\}\cup (P_x(\alpha, \beta, \phiv)-y)$
where $P_x(\alpha, \beta, \phiv)$
 is the $(\alpha, \beta)$-chain containing both $x$ and $y$.
\item {\bf Adding edges by pairs:} Repeatedly pick two boundary edges $f_1$ and $f_2$ of $T^*$ with  $\phiv(f_1) = \phiv(f_2) \in \phibar(V(T^*))$ and redefine  $T^* := T^*\cup \{f_1, f_2\}$ until $T^*$ is closed. \end{itemize}

The path $P_x(\alpha, \beta, \phiv)$ in the above definition is called the {\it trunk} of $T^*$ and $h(T^*):=|V(P_x(\alpha, \beta, \phiv))|$ is called the {\it height} of $T^*$.

\begin{LEM}\label{LEM:balanced-T} {\em [Scheide~\cite{Scheide-2010}]}
Let $G$ be a $k$-critical graph with $k \ge \Delta +1$ and $T\in \mathcal{T}$ be a balanced Tashkinov tree of a $k$-triple $(G, e, \phiv)$ with $h(T)$ being maximum. Then, $h(T)\ge 3$ is odd. Moreover, if $h(T) =3$ then $G$ is elementary.
\end{LEM}

\begin{COR}\label{COR:T-order-exit-vertex}
Let $G$ be a non-elementary $k$-critical graph with $k\ge \Delta +1$  and $T\in \mathcal{T}$ be a balanced Tashkinov tree of a $k$-triple $(G, e, \phiv)$ with $h(T)$ being maximum. Then
$|T|\ge 2(k-\Delta)+1$.
\end{COR}
\pf Since $G$ is not elementary, $T$ is not strongly closed with respect to $\phiv$. There is an exit vertex $v$
by Lemma~\ref{LEM:PassAll}, so $\phibar(v) \subseteq \phiv(E(T))$ by
 Lemma~\ref{LEM:exit-color}. Since $T$ is balanced and $h(T) \ge 5$ by Lemma \ref{LEM:balanced-T}, each used color is assigned to at least two edges of $E(T)$.  Thus,
 \[
 |T| = ||T|| +1 \ge 2 |\phibar(v)| +1 \ge 2(k-\Delta) +1. \qed
 \]

Working on balanced Tashkinov trees, Scheide proved the following result.
\begin{LEM}\label{LEM:Small-T}{\em [Scheide ~\cite{Scheide-2010}]}
Let $G$ be a $k$-critical graph with $k \ge \Delta +1$. If $|T| <11$ for all Tashkinov trees $T$, then $G$ is  elementary.
\end{LEM}

\section{An extension of  Tashkinov trees}
\subsection{Definitions and basic properties}
In this section, we always assume that $G$ is a {\it non-elementary} $k$-critical graph with $k\ge \Delta+1$ and $T_0\in \mathcal{T}$ is a maximum Tashkinov tree of $G$. Moreover, we assume that $T_0$ is a Tashkinov tree of the $k$-triple $(G, e, \phiv)$.
\begin{DEF}\label{DEF:stable}
Let  $\varphi_1,\varphi_2\in \mathcal{C}^k(G-e)$ and  $H\subseteq G$ such that $e\in E(H)$. We say that $H$  is {\it $(\varphi_1,\varphi_2)$-stable} if  $\varphi_1(f) = \varphi_2(f)$ for every $f\in E(G[V(H)])\cup \partial(H)$, that is, $\varphi_1(f)\neq \varphi_2(f)$ implies that $f\in E(G-V(H))$.
\end{DEF}
Following the definition, if a Tashkinov tree $T_0$ of $(G,e,\varphi_1)$ is $(\varphi_1,\varphi_2)$-stable, then it is also a Tashkinov tree of $(G, e,\varphi_2)$. Moreover, the sets of missing colors of $T_0$, used colors of $T_0$, and free colors of $T_0$ are the same in both colorings $\varphi_1$ and $\varphi_2$.

The following definition of {\it connecting edges} will play a critical role in our extension based on a maximum Tashkinov tree.
\begin{DEF}\label{DEF:Conn}
Let $H\subseteq G$ be a subgraph such that $T_0\subseteq H$. A color $\delta$ is called a {\it defective color} of $H$ if $H$ is closed, $\delta\not\in \phibar(V(H))$ and $|E_{\delta}(\partial(H))| \ge 2$.
Moreover, an edge $f\in \partial(H)$ is called a {\it connecting edge} if $\delta :=\phiv(f)$ is a defective color of $H$  and there is a missing color $\gamma\in \phibar(V(T_0)) -\phiv(E(H))$ of  $T_0$ such that the following two properties hold.
\begin{itemize}
\item The $(\gamma, \delta)$-chain $P_u(\delta, \gamma, \phiv)$ contains all edges in $E_{\delta}(\partial(H))$, where $u$ is the unique vertex in $V(T_0)$ such that $\gamma \in \phibar (u)$;
\item Along the linear order $\preceq_{(u, P_u(\gamma, \delta, \phiv))}$, $f$ is the first boundary edge on $P_u(\gamma, \delta, \phiv)$ with color $\delta$.
\end{itemize}
\end{DEF}

In the above definition, we call the successor $f^s$ of $f$  along $\preceq_{(u, P_u(\gamma, \delta, \phiv))}$ the {\it companion}  of $f$, $(f,f^s)$ a {\it connecting edge pair} and $(\delta,\gamma)$ a {\it connecting color pair}.  Since $P_u(\gamma, \delta, \phiv)$ contains all edges in $E_{\delta}(\partial(H))$, we have that $f^s$ is not incident to any vertex in $H$ and $\phiv(f^s) = \gamma$.

\begin{DEF}\label{DEF:ETT}
 We call a tree $T$ an {\bf Extension of a Tashkinov Tree (ETT)} of $(G, e, \phiv)$ based on $T_0$ if $T$ is incrementally obtained from $T:=T_0$ by  repeatedly adding edges to  $T$ according to the following two operations subject to $\Gamma^f(T_0) - \phiv(E(T)) \ne \emptyset$:
 \begin{itemize}
 \item {\bf ET0:} If $T$ is closed,  add a connecting edge pair $(f, f^s)$, where $\varphi(f)$ is a defective color and $\varphi(f^s)\in \Gamma^f(T_0) - \phiv(E(T))$,  and rename $T:=T\cup \{f,f^{s}\}$.
 \item {\bf ET1:}  Otherwise, add an edge $f\in \partial(T)$ with $\phiv(f) \in \phibar(V(T))$ being a missing color of $T$, and rename $T:=T\cup \{f\}$.
\end {itemize}
 \end{DEF}

Note that the above extension algorithm ends with  $\Gamma^f(T_0)\subseteq \varphi(E(T))$.
Let $T$  be an ETT of   $(G, e, \phiv)$. Since $T$ is defined incrementally from $T_0$,  the edges added to $T$ follow a linear order $\prec_{\ell}$.  Along the linear order $\prec_{\ell}$, for any initial subsequence $S$ of $E(T)$, $T_0\cup S$ induces a tree; we call it a {\it premier segment} of $T$ provided that when a connecting edge is in $S$, its companion must be in $S$. Let $f_1, \, f_2, \dots, f_{m+1}$ be all connecting edges  with  $f_1\prec_{\ell} f_2\prec_{\ell} \dots \prec_{\ell} f_{m+1}$.  For each $1\le i\le m+1$,  let $T_{i-1}$ be the premier subtree induced by $T_0$ and edges before $f_i$ in the ordering $\prec_{\ell}$.    Clearly, we have  $T_0\subset T_1 \subset T_2 \subset \dots \subset T_m \subset T$. We call $T_i$ a {\it closed segment} of $T$ for each $0\le i\le m$, $T_0\subset T_1 \subset T_2 \subset \dots \subset T_m \subset T$ the {\it ladder of $T$}, and $T$ an {\it ETT with $m$-rungs}. We use $m(T)$ to denote the number of rungs of $T$.  For each edge $f\in E(T)$ with $f\ne e$,  following the linear order $\prec_{\ell}$,  the end of $f$ is called the {\it in-end} if it is in $T$ before $f$ and the other one is called the   {\it out-end} of $f$. For any edge $f\in E(T)$, the subtree induced by $T_0$, $f$ and all its predecessors  is called an $f$-{\it segment} and denoted by $Tf$.

 Let $\mathbb{T}$ denote the set of all ETTs based on $T_0$.  We now define a binary relation $\prec_t$  of $\mathbb{T}$ such that for  two $T, T^*\in \mathbb{T}$, we call $T\prec_t T^*$ if either $T=T^*$ or there exists $s$ with $1\le s\le \min\{m+1,m^*+1\}$ such that $T_h = T^*_h$ for every $0 \le h< s$ and $T_s \subsetneq T^*_s$, where $T_0 \subset T_1 \subset \dots \subset T_s \subset \dots \subset T_{m} \subset T_{m+1}(=T)$ and $T_0^*(=T_0) \subset T^*_1\subset \dots \subset T^*_s \subset \dots \subset T^*_{m^*+1} (=T^*)$ are the ladders of $T$ and $T^*$, respectively. Notice that in this definition, we only consider the relations of $T_h$ and $T^*_h$ for $h\leq s$. Clearly, for any three ETTs $T$, $T'$ and $T^*$, $T\prec_t T'$ and $T' \prec_t T^*$ give $T\prec_t T^*$. So, $\mathbb{T}$ together with $\prec_t$ forms a poset, which is denoted by $(\mathbb{T}, \prec_t)$.

 \begin{LEM}\label{LEM:Max}
 In the poset $(\mathbb{T}, \prec_t)$, if $T$ is a maximal tree over all ETTs with at most $|T|$ vertices, then any premier segment $T'$ of $T$  is also a maximal tree over all ETTs with at most $|T'|$ vertices.
 \end{LEM}
\proof Suppose on the contrary: {\it there is a premier segment $T'$ of $T$ and an ETT $T^*$ with $|T^*| \le |T'|$ and $T'\prec_t T^*$}.  We assume that $T' \ne T^*$.  Let $T_0\subset T_1\subset \dots \subset T_{m'}\subset T'$  and $T_0 \subset T_1^* \subset \dots \subset T_{m^*}^*\subset T^*$ be the ladders of $T'$ and $T^*$, respectively.  Since $T'\prec_t T^*$, there exists $s$ with $1\le s\le \min\{m'+1, m^*+1\}$ such that
$T_j = T^*_j$ for each $0\le j\le s-1$ and $T_s \subsetneq T^*_s$, where $T'_{m'+1}=T'$ and $T_{m^*+1}^* = T^*$.  Since $|T^*| \le |T'|$, we have $s <m'+1$.  Since $T'$ is a premier segment of $T$,  $T_0\subset T_1\subset \dots \subset T_{m'}$ is a part of the ladder of $T$.  So, we have $T\prec_t T^*$, giving a contradiction to the maximality of $T$. \qed

\begin{LEM}\label{LEM:FixConn}
Let $T$ be a maximal ETT in $(\mathbb{T}, \prec_t)$ over all ETTs with at most $|T|$ vertices, and let $T_0 \subset T_1 \subset \cdots \subset T_m \subset T$ be the ladder of $T$. Suppose $T$ is an ETT of $(G,e,\phiv_1)$. Then for every $\phiv_2 \in \mathcal{C}^k(G-e)$ such that $T_m$ is $(\varphi_1,\varphi_2)$-stable, $T_m$ is an ETT of  $(G, e, \varphi_2)$. Furthermore, if $T_m$ is elementary, then for every $\gamma \in \Gamma^f(T_0) - \varphi_1(E(T_m))$ and $\delta \not \in \phibar_1(V(T_m))$, $P_u(\gamma, \delta,\varphi_2) \supseteq \partial_{\delta}(T_m)$ where $u \in V(T_0)$ such that $\gamma \in \phibar_1(u)$.
\end{LEM}
\proof
Suppose on the contrary: let $T$  be a counterexample to Lemma~\ref{LEM:FixConn} with minimum number of vertices. Let $T_0 \subset \cdots \subset T_m \subset T$ be the ladder of $T$ and let $\phiv_1,\phiv_2 \in \mathcal{C}^k(G-e)$ be two edge colorings such that $T$ is an ETT of $(G,e,\phiv_1)$, $T_m$ is $(\phiv_1,\phiv_2)$-stable and either
\begin{itemize}
\item[(1)] $T_m$ is not an ETT of $(G,e,\phiv_2)$ or
\item[(2)] $T_m$ is elementary and there exist $\gamma \in \Gamma^f(T_0) - \varphi_1(E(T_m))$ and $\delta \not \in \phibar_1(V(T_m))$ such that $P_u(\gamma, \delta,\varphi_2) \not\supseteq \partial_{\delta}(T_m)$ where $u \in V(T_0)$ such that $\gamma \in \phibar_1(u)$.
\end{itemize}

By the minimality of $T$, we observe that $|T|=|T_m|+2$. Furthermore, since $T_0 \in \mathcal{T}$ is a maximum Tashkinov tree of $G$, it follows that $m\ge 1$ by Lemma~\ref{LEM:PassAll}.

First, we show that (1) does not hold, in other words, $T_m$ is an ETT of $(G,e,\varphi_2)$.  Since colors for edges incident to vertices in $T_m$ are the same in both $\varphi_1$ and $\varphi_2$, we only need to show that each connecting edge pair  in coloring $\varphi_1$ is still a connecting edge pair in coloring $\varphi_2$. For $0\le j \le m-1$ let $(f_j,f_j^s)$ be the connecting edge pair of $T_j$ and let $(\delta_j,\gamma_j)$ be the corresponding connecting color pair with respect to $\varphi_1$.
Since $T_{j+1}$ is $(\varphi_1,\varphi_2)$-stable and an ETT of $(G,e,\varphi_1)$ and $T_{j+1} \subsetneq T$, by the minimality of $T$, it follows that
$P_{u_j}(\gamma_j,\delta_j,\varphi_2)$ contains $\partial_{\delta_j}(T_j)$ where $u_j$ is the unique vertex in $V(T_0)$ with $\gamma_j \in \phibar_1(u_j)$.
Moreover, since $T_{j+1}$ is $(\varphi_1,\varphi_2)$-stable, it follows that $f_j$ is the first boundary edge on $P_{u_j}(\gamma_j,\delta_j,\varphi_2)$ with color $\delta_j$ and $f_j^s$ being its companion.
So $(f_j,f_j^s)$ is still a connecting edge pair  in $\varphi_2$. We point out that $P_{u_j}(\gamma_j,\delta_j,\varphi_1)$ and $P_{u_j}(\gamma_j, \delta_j, \varphi_2)$ may be different in $(G, e, \varphi_1)$ and $(G, e, \varphi_2)$.

Thus (2) holds and there exist $\gamma \in \Gamma^f(T_0) - \varphi_1(E(T_m))$ and $\delta \not \in\phibar_1(V(T_m))$ such that $P_u(\gamma, \delta,\varphi_2) \not\supseteq \partial_{\delta}(T_m)$.
Let $P=P_u(\gamma, \delta,\varphi_2)$.
Since $T_m$ is both elementary and closed and $u$ is one of the two ends of $P$, the other end of $P$ must be in $V \setminus V (T_m)$.
So, $E(P)\cap E_{\delta}(\partial(T_m)) \neq \emptyset$. Let $Q$ be another $(\gamma, \delta)$-chain such that $E(Q) \cap E_{\delta}(\partial(T_m))\neq  \emptyset$. Let $\varphi_3 := \varphi_2/Q$ be a coloring of $G- e$ obtained from $\varphi_2$ by interchanging colors assigned on $E(Q)$.

Let $(f,f^s)$ be the connecting edge pair of $T_{m-1}$, and $T'=T_{m-1} \cup \{f,f^s\}$. We claim that $E(T') \cap E(Q) =\emptyset$.
By the minimality of $T$, $P$ contains every edge of $E_{\delta}(\partial(T_{m-1}))$, and so $E(T_{m-1}) \cap E(Q)=\emptyset$.  If $\varphi_2(f) \neq \delta$ then $f \not \in E(Q)$ and if $\varphi_2(f)=\delta$ then $f \in E(P)$ so $f\not \in E(Q)$. Thus $f \not \in E(Q)$.
Lastly, $\varphi_2(f^s) \neq \delta$ since $\delta \in \phibar_2(V(T_m))$ and $\varphi_2(f^s) \neq \gamma$ since $\gamma \not \in \varphi_2(E(T_m))$, so $f^s \not \in E(Q)$.

Observe that $T'$ is an ETT of $(G,e,\varphi_1)$ with ladder $T_0 \subset \cdots \subset T_{m-1}$ and is $(\varphi_1,\varphi_3)$-stable. Moreover $|T'| \le |T_m| <|T|$. Therefore, by the minimality of $T$, $T_{m-1}$ is an ETT of $(G,e,\varphi_3)$, and  because we do not use any edge in $Q$ when we extend $T_{m-1}$ to $T_m$, $T_m$ is also an ETT of $(G,e,\varphi_3)$ which is not closed. However, it is a contradiction that $T$ is a maximal ETT. \qed

In Lemma \ref{LEM:FixConn}, by taking $\phiv_1=\phiv_2$, we easily obtain the following lemma.



 \begin{LEM} \label{LEM:EXconn}
Let $T$ be a maximal ETT in $(\mathbb{T}, \prec_t)$ over all ETTs with at most $|T|$ vertices, and let $T_0\subset T_1\subset \dots \subset T_m\subset T$ be the ladder of $T$. Suppose $T$ is an ETT of $(G,e,\varphi)$.  If $T_m$ is elementary and $\Gamma^f(T_0) -\phiv(E(T))\ne \emptyset$,
then for any $\gamma\in \Gamma^f(T_0) -\phiv(E(T))$ and $\delta\notin \phibar(V(T_m))$, $P_u(\gamma, \delta, \phiv) \supset E_{\delta}(\partial (T_i))$ for every $i$ with $0\le i\le m$, where $u\in V(T_0)$ such that $\gamma\in \phibar(u)$.
\end{LEM}


\begin{LEM}\label{LEM:ETT>F}
For every  ETT $T$ of $(G,e,\varphi)$ based on $T_0$, if $T$ is elementary such that $|\Gamma^f(T_0)| >m(T)$ and  $ |E(T) - E(T_0)|-m(T) < |\phibar(V(T_0))|$, then there exists an ETT $T^*$ containing $T$ as a premier segment.
\end{LEM}
\proof Let $T$ be an ETT of $(G,e,\varphi)$ and $m = m(T)$.  Since $\phiv(f_i) \notin \phibar(V(T_0))$ for each connecting edge $f_i$, where $i\in \{1,2,\cdots,m\}$,  we have $|\phiv(E(T)-E(T_0))\cap \phibar(V(T_0))| \le |E(T) - E(T_0)| -m < |\phibar(V(T_0))|$. So,  $\phibar(V(T_0)) - \phiv(E(T)-E(T_0)) \ne \emptyset$.  Let $\gamma\in\phibar(V(T_0)) - \phiv(E(T)-E(T_0))$.

We may assume $\gamma \notin \phiv(E(T_0))$, i.e., $\gamma \in  \Gamma^f(T_0)$. Since $m < |\Gamma^f(T_0)|$, there exists a color $\beta\in \Gamma^f(T_0) -\{\gamma_1, \gamma_2, \dots, \gamma_m\}$.  Since $T_0$ is closed, a $(\beta, \gamma)$-chain is either in $G[V(T_0)]$ or vertex disjoint from $T_0$. Let $\varphi_1$ be obtained from $\phiv$ by interchanging $\beta$ and $\gamma$ for edges in $E_{\beta}(G-V(T_0))\cup E_{\gamma}(G-V(T_0))$.   Clearly, $T_0$ is $(\phiv, \varphi_1)$-stable. So,
$T$ is also an ETT of $(G, e, \phiv_1)$.
Since $\gamma\notin \phiv(E(T)-E(T_0))$, we have $\beta\notin \phiv_1(E(T))$, so the claim holds.

We can apply {\bf ET0} and {\bf ET1} to extend $T$ to a larger tree $T^*$ unless $T$ is closed and does not have a connecting edge. In this case, $T$ is both elementary and closed. Since $G$ itself is not elementary, $T$ is not strongly closed. Thus, $T$ has a defective color $\delta$.
Since $T$ does not have a connecting edge, $P_v(\gamma, \delta, \phiv)$ does not contain all edges of $E_{\delta}(\partial(T))$, where $v\in V(T_0)$ is the unique vertex with $\gamma\in \phibar(v)$. Let $Q$ be another $(\gamma, \delta)$-chain containing some edges in $E_{\delta}(\partial(T))$ and let $\varphi_2= \phiv/Q$.  By Lemma~\ref{LEM:EXconn}, $Q$ is disjoint from $T_m$, where $T_m$ is the largest closed segment of $T$. So, $T_m$ is $(\phiv, \varphi_2)$-stable.  By Lemma~\ref{LEM:FixConn}, $T_m$ is an ETT of $(G, e, \varphi_2)$, which in turn gives that $T$ is also an ETT of $(G, e, \varphi_2)$.  Applying {\bf ET1},  we extend $T$ to a larger ETT $T^*$, which contains $T$ as a premier segment.
\qed

\subsection{The major result}
The following result is fundamental for both Theorems~\ref{THM:cubic} and \ref{THM:Jm19}.
\begin{THM}\label{LEM:Elem}
Let $G$ be a $k$-critical graph with $k\ge \Delta +1$ and $T$ be a maximal  ETT  over all ETTs with at most $|T|$ vertices in the poset $(\mathbb{T}, \prec_t)$. Suppose $T$ is an ETT of $(G,e,\varphi)$. If $|E(T) - E(T_0)| -m(T) < |\phibar(V(T_0))|-1$ and $m(T) < |\Gamma^f(T_0)|-1$, then $T$ is elementary.
\end{THM}
\proof
Suppose on the contrary: {\it let $T$ be a counterexample to Theorem~\ref{LEM:Elem} with minimum number of vertices.} And we assume that  $(G, e, \phiv)$ is the triple in which $T$ is an ETT.

By Theorem~\ref{THM:TashOrigi}, we have  $T\supsetneq T_0$.
For any premier segment $T'$ of $T$, by Lemma~\ref{LEM:Max},    $T'$ is maximal over all ETTs with at most $|T'|$ vertices.   Additionally, following the definition,  we can verify that $|E(T') -E(T_0)| -m(T') \le |E(T) - E(T_0)| -m(T)$ and $m(T')\leq m(T)$. So, every premier segment of $T$ satisfies the conditions of Theorem~\ref{LEM:Elem}.  Hence, Theorem~\ref{LEM:Elem} holds for all premier segments of $T$ which are proper subtrees of $T$.   Let $T_0\subset T_1\subset \dots \subset T_m\subset T$ be the ladder of $T$.

Let $v_1, v_2$ be two distinct vertices in $T$ such that there is a color $\alpha \in \phibar(v_1)\cap \phibar(v_2)$.  For each connecting edge $f_i$ with $1\le i \le m$,  let $(\delta_i, \gamma_{\delta_i})$ denote the corresponding color pair, where $\phiv(f_i) = \delta_i$.  According to the definition of ETT, $\gamma_{\delta_1}, \gamma_{\delta_2}, \dots, \gamma_{\delta_m}$ are pairwise distinct while $\delta_1$, $\delta_2$, $\dots$, $\delta_m$ may not be.  Let $L=\{\gamma_{\delta_1}, \gamma_{\delta_2}, \dots, \gamma_{\delta_m}\}$. In the paper \cite{CYZ-2011} by Chen et al., the condition $\phibar(v)\not\subseteq L$ is needed for any $v\in V(T)-V(T_0)$. In the following proof, we overcome this constraint. We make the following assumption.

{\flushleft \bf Assumption 1:}  We assume that over all colorings in $\mathcal{C}^k (G-e)$ such that $T$ is a minimum counterexample, the coloring $\phiv\in \mathcal{C}^k(G-e)$ is one such that $|\phibar(V(T_0)) -(\phiv(E(T)-E(T_0))\cup\{\alpha\})|$ is minimum.

The following claim states that we can use other missing colors of $T_0$ before using free colors of $T_0$ except those in $L$.

\begin{CLA}\label{CLA:usedfirst}
We may assume that if $\phiv(E(T)-E(T_0))\cap (\Gamma^f(T_0) -(L\cup\{\alpha\})) \ne \emptyset$, then $\phiv(E(T)-E(T_0)) \supset \phibar(V(T_0)) -\Gamma^f(T_0)$.  \end{CLA}

\proof  Assume that there is a color  $\gamma \in \phiv(E(T)-E(T_0))\cap (\Gamma^f(T_0) -(L\cup\{\alpha\}))$ and there is a color $\beta\in (\phibar(V(T_0)) -\Gamma^f(T_0)) - \phiv(E(T)-E(T_0))$. Since $T_0$ is closed, a $(\beta, \gamma)$-chain is either in $G[V(T_0)]$ or disjoint from $V(T_0)$. Let $\varphi_1$ be obtained from $\phiv$ by interchanging colors $\beta$ and $\gamma$ on all $(\beta, \gamma)$-chains disjoint from $V(T_0)$.  It is readily seen that $T_0$ is $(\phiv, \varphi_1)$-stable. Since both $\gamma$ and $\beta$ are in $\phibar(V(T_0))-L$,  $T$ is also an ETT of $(G, e, \varphi_1)$.  In coloring $\varphi_1$, we still have $\gamma\in\Gamma^{f}(T_0)-(L\cup \{\alpha\})$ and $\beta\in \phibar_1(V(T_0))-\Gamma^{f}(T_0)$. However, $\gamma$ is not used on $T-T_0$ while $\beta$ is used. Additionally,  Assumption 1 holds since  $|\phibar(V(T_0)) -(\phiv(E(T)-E(T_0))\cup\{\alpha\})|  =  |\phibar_1(V(T_0)) -(\phiv_1(E(T)-E(T_0))\cup\{\alpha\})|$.
 By repeatedly applying this argument, we show that  Claim~\ref{CLA:usedfirst} holds. \qed

Since $m(T) < |\Gamma^f(T_0)|-1$,  we have $\Gamma^f(T_0) - (L\cup\{\alpha\}) \ne \emptyset$.  Since $|E(T)-E(T_0)| -m(T) < |\phibar(V(T_0))|-1$, we have  $\phibar(V(T_0)) -(\phiv(E(T)-E(T_0))\cup\{\alpha\}) \ne \emptyset$.  By Claim~\ref{CLA:usedfirst}, we have the following claim.
\begin{CLA}\label{CLA:Non0}
We may assume that $\Gamma^f(T_0) - (\phiv(E(T))\cup \{\alpha\}) \ne \emptyset$.
\end{CLA}

We consider two cases to complete the proof according to the type of the last operation in adding edge(s)
to extend $T_0$ to $T$.

{\flushleft \bf Case 1:}  The last operation is {\bf ET0}, i.e., the two edges in the connecting edge pair $(f, f^s)$ are the last two edges in $T$ following  the linear order $\prec_{\ell}$.

Let $x$ be the in-end of $f$, $y$ be the out-end of $f$ (in-end of $f^s$), and $z$ be the out-end of $f^s$. In this case, we have $V(T)=V(T_m)\cup\{y, z\}$, i.e., $T'=T_m$.   Let $\delta = \phiv(f)$ be the defective color and $\gamma_{\delta}\in \Gamma^f(T_0) -\phiv(E(T_m))$ such that $f$ is the first edge in $\partial(E(T_m))$ along $P:=P_u(\gamma_{\delta}, \delta, \phiv)$ with color $\delta$, where $u\in V(T_0)$ such that $\gamma_{\delta}\in \phibar(u)$. Recall that $v_1$ and $v_2$ are the two vertices in $T$ such that $\alpha \in \phibar(v_1)\cap\phibar(v_2)$. We have $\{v_1, v_2\}\cap \{y, z\} \ne \emptyset$. We consider the following three subcases  to lead a contradiction.

{\bf \noindent Subcase 1.1: $\{v_1, v_2\} = \{y, z\}$ }.

Assume, without loss of generality, $y=v_1$ and $z=v_2$.
Since $f^s$ is the successor of $f$ along the linear order $\preceq_{(u,P)}$, $\phiv(f^s) = \gamma_{\delta}$.  So, $f^s$ is an $(\alpha, \gamma_{\delta})$-chain. Let $\varphi_1 = \phiv/f^s$, a coloring obtained from $\phiv$ by changing color on $f^s$ from $\gamma_{\delta}$ to $\alpha$. Then $T_m$ is $(\varphi,\varphi_1)$-stable. By Lemma \ref{LEM:FixConn}, $T_m$ is an ETT of $(G,e,\varphi_1)$ and $\gamma_{\delta}$ is missing at $y$ in $\varphi_1$, which in turn gives that $P_u(\gamma_{\delta},\delta,\varphi_1):=uPy$ only contains one edge $f\in E_{\delta}(\partial(T_m))$, giving a contradiction to Lemma \ref{LEM:EXconn}.

{\bf \noindent Subcase 1.2: $\alpha\in (\phibar(y)-\phibar(z))\cap \phibar(V(T_{m}))$.}

Since $\delta, \gamma_{\delta} \in \phiv(y)$ and $\alpha\in \phibar(y)$,   $\alpha\notin\{\delta, \gamma_{\delta}\}$.
We may assume that $\alpha\in \Gamma^f(T_0)-\phiv(E(T))$.
Otherwise,  let $\beta\in\Gamma^f(T_0)-\phiv(E(T))$ and consider the $(\alpha,\beta)$-chain $P_1:=P_y(\alpha,\beta,\phiv)$.
Since $\alpha,\beta\in \phibar(V(T_{m}))$
and $V(T_{m})$ is closed with respect to $\phiv$ by the assumption,
we have $V(P_1)\cap V(T_{m})=\emptyset$. Let $\phiv_1=\phiv/P_1$.
Since  $\{\alpha,\beta\}\cap\{\delta, \gamma_{\delta}\}=\emptyset$, we have
$f^s\notin E(P_1)$. Hence $T_{m}$ is $(\phiv, \phiv_1)$-stable, which gives that
$T_m$ is an ETT of $(G, e,\phiv_1)$, so is $T$.  The claim follows from  $\beta\in\phibar_1(y)
\cap (\Gamma^f(T_0)-\phiv_1(E(T)))$.

Consider the $(\alpha,\gamma_{\delta})$-chain $P_2:=P_y(\alpha,
\gamma_{\delta},\phiv)$.  Since $\alpha, \gamma_{\delta}\in \phibar(V(T_0))$ and $T_m$ is closed,
$V(P_2)\cap V(T_{m})=\emptyset$. Let $\phiv_2=\phiv/P_2$.  Clearly, $T_m$ is $(\phiv, \phiv_2)$-stable. By Lemma~\ref{LEM:FixConn}, $T_m$ is an ETT of $(G, e, \phiv_2)$, so is $T$. Then
$P_{u}(\gamma_{\delta},\delta,\phiv_2)$ is the subpath of  $P_u(\gamma_{\delta}, \delta, \phiv)$ from $u$ to $y$. So, it does not contain
all edges in
$E_{\delta}(\partial(T_{m}))$, which gives a contradiction to Lemma~\ref{LEM:EXconn}.

{\bf \noindent Subcase 1.3: $\alpha\in
(\phibar(z)-\phibar(y))\cap \phibar(V(T_{m}))$.}

Since $P_u(\gamma_{\delta}, \delta, \phiv)$
 contains all the edges in
$E_{\delta}(\partial(T_{m}))$ and $\alpha\in \phibar(z)$,  we have
$\alpha\notin\{\delta, \gamma_{\delta}\}$.
Following a similar argument given in Subcase 1.2, we may assume that
$\alpha\in\Gamma^{f}(T_0)-\phiv(E(T))$.  Let $v$ be the unique vertex in $V(T_0)$ with $\alpha\in \phibar(v)$. Let $\beta\in \phibar(y)$, $P_v:=P_{v}(\alpha,\beta,\phiv)$, $P_y:=P_{y}(\alpha,\beta,\phiv)$ and $P_z:=P_{z}(\alpha,\beta,\phiv)$. We claim that $P_v=P_y$. Suppose, on the contrary, that $P_v\neq P_y$. By Lemma \ref{LEM:EXconn}, $E(P_v)\supset E_{\beta}(\partial(T_{m}))$. Therefore, $V(P_y)\cap V(T_{m})=\emptyset$. Let $\phiv_1=\phiv/P_y$. In $(G, e, \phiv_1)$,  $T$ is an ETT and $\alpha\in \phibar_1(y)\cap\phibar_1(V(T_{0}))$. This leads back to either Subcase 1.1 or Subcase 1.2. Hence, $P_v=P_y$ and it is vertex disjoint with $P_z$. Let $\phiv_2=\phiv/P_z$.  By Lemma \ref{LEM:EXconn}, $E(P_v)\supset E_{\beta}(\partial(T_m))$.   So, $V(P_z)\cap V(T_m)=\emptyset$, which in turn gives that $T$ is an ETT of $(G, e, \phiv_2)$ and $\beta\in \phibar_2(y)\cap\phibar_2(z)$. This leads back to Subcase 1.1.

{\flushleft \bf Case 2:} The last edge $f$ is added to $T$ by {\bf ET1}.

Let $y$ and $z$ be the in-end and out-end of $f$, respectively, and let $T' =T-z$. Clearly, $T'$ is a premier segment of $T$ and $T_m \subsetneq T'$. In this case, we assume that $z=v_2$, i.e., $\alpha\in  \phibar(z)\cap \phibar(v_1)$ and $v_1\in V(T')$. Recall that $v_1$ and $v_2$ are the two vertices in $T$ such that $\alpha\in\phibar(v_1)\cap\phibar(v_2)$.

\begin{CLA}\label{CLA:OneComp}
For any color $\gamma\in \Gamma^f(T_0)$ and any color $\beta\in \phibar(V(T'))$, let $u\in V(T_0)$ such that
$\gamma\in \phibar(u)$ and $v\in V(T')$ such that $\beta\in \phibar(v)$. Denote by $e_v\in E(T)$ the edge containing $v$ as the out-end and $e_v\prec_{\ell} e^*$ for every $e^*\in E(T)$ with $\varphi(e^*)=\gamma$, then $u$ and $v$ are on the same $(\beta, \gamma)$-chain.
\end{CLA}
\proof Since $T_m$ is both elementary and closed, $u$ and $v$ are on the same $(\beta, \gamma)$-chain if $v\in V(T_m)$.  Suppose $v\in V(T) -V(T_m)$ and, on the contrary, $P_u:=P_u(\gamma, \beta, \phiv)$ and $P_v:=P_v(\gamma, \beta, \phiv)$ are vertex disjoint. By Lemma~\ref{LEM:EXconn}, $E(P_u)\supset E_{\beta}(\partial(T_m))$, so $V(P_v)\cap V(T_m) = \emptyset$.  Let $\phiv_1 = \phiv/P_v$ be the coloring obtained by interchanging the colors $\beta$ and $\gamma$ on $P_v(\gamma, \beta, \phiv)$. Clearly, $T_m$ is $(\phiv, \phiv_1)$-stable. By Lemma~\ref{LEM:FixConn}, $T_m$ is an ETT of $(G, e, \phiv_1)$. As $e_v\prec_{\ell} e^*$ for every $e^*\in E(T)$ with $\varphi(e^*)=\gamma$, we can extend $T_m$ to $Te_{v}$ such that $Te_{v}$ is still an ETT of $(G,e,\varphi_1)$.  But, in the coloring $\phiv_1$, $\gamma\in \phibar_1(u)\cap \phibar_1(v)$, which gives a contradiction to the minimality of $|T|$. \qed

\begin{CLA}\label{CLA:a-free}
We may assume $\alpha\in \Gamma^f(T_0)-\phiv(E(T_{m}))$.
\end{CLA}
\proof Otherwise, by Claim~\ref{CLA:Non0}, let $\gamma\in \Gamma^f(T_0) -(\phiv(E(T))\cup\{\alpha\})$.
Let $\phiv_1$ be obtained from $\phiv$ by interchanging colors $\alpha$ and $\gamma$ for edges in $E_{\alpha}(G-V(T_m))\cup E_{\gamma}(G-V(T_m))$.   Since $T_m$ is closed, $\phiv_1$ exists.  Clearly, $T_m$ is $(\varphi,\varphi_1)$-stable. By Lemma~\ref{LEM:FixConn}, $T_m$ is an ETT of $(G, e, \phiv_1)$, so is $T$. In the coloring $\phiv_1$, $\gamma\in \phibar_1(z)$ but is not used on $T_m$.
\qed

Applying Claim~\ref{CLA:Non0} again if it is necessary,  we assume both Claim~\ref{CLA:Non0} and Claim~\ref{CLA:a-free} hold.  Recall that $z$ is the out-end of $f$ and $y$ is the in-end of $f$, and $\alpha \in \phibar(v_1)\cap \phibar(z)$.

{\bf \noindent Subcase 2.1:} $y\in V(T') - V(T_m)$, i.e., $f\notin \partial(T_m)$.

\begin{CLA}\label{CLA:a-used}
Color $\alpha$ is used in $E(T -T_m)$, i.e., $\alpha\in \phiv(E(T-T_m))$.
\end{CLA}
\proof Suppose on the contrary that $\alpha \notin \phiv(E(T-T_m))$.  By Claim~\ref{CLA:a-free}, we may assume that $\alpha \notin \phiv(E(T_m))$,   so $\alpha \notin \phiv(E(T))$. Let $\phiv(f)=\theta$ and  $\beta \in \phibar(y)$ be a missing color of $y$.    We consider the following two cases according to whether $y$ is the last vertex of $T'=T-z$.

We first assume that $y$ is the last vertex of $T'$.   Let
$P_{v_1}:=P_{v_1}(\alpha,\beta,\phiv)$, $P_y:=P_y(\alpha,\beta,\phiv)$ and $P_z:=P_z(\alpha,\beta,\phiv)$ be
 $(\alpha, \beta)$-chains containing vertices $v_1$, $y$ and $z$, respectively.
  By Claim~\ref{CLA:OneComp},  we have $P_{v_1}=P_y$, so it is disjoint from $P_z$.
By Lemma~\ref{LEM:EXconn}, $E(P_{v_1})\supset E_{\beta}(\partial(T_m))$, so $V(P_z)\cap V(T_m) = \emptyset$. Let $\phiv_1 = \phiv/P_z$ be the coloring obtained from $\phiv$ by interchanging colors $\alpha$ and $\beta$ on $P_z$.  Since $\alpha\notin \phiv(E(T-T_m))$ and $\beta\in \phibar(y)-\phibar(V(T'))$, $\beta\not\in \varphi_1(E(T-T_m))$.  Clearly, $T_m$ is $(\phiv, \phiv_1)$-stable. By Lemma~\ref{LEM:FixConn}, $T_m$ is an ETT of $(G, e, \phiv_1)$, so is $T$.  In the coloring $\phiv_1$, $\theta = \phiv_1(f)$ and $f$ itself is a $(\beta, \theta)$-chain. Let $\phiv_2 = \phiv_1/f$ be the coloring obtained from $\phiv_1$ by changing color $\theta$ to $\beta$ on $f$.  Since $f$ is disjoint from $T_m$,  we can verify that $T$ is an ETT of $(G, e, \phiv_2)$ by applying Lemma~\ref{LEM:FixConn}. Since $f$ is not a connecting edge, $\theta\in \phibar(V(T'))$, which in turn shows that $T'$ is not elementary with respect to $\varphi_2$, giving a contradiction to the minimality of $|T|$.

We now assume that $y$ is not the last vertex of $T'$;  and let $x$ be the last one.  Recall $\theta = \phiv(f)$.  If $\theta\in\varphi(x)$ then $T-x$ is not an elementary ETT of $(G,e,\varphi)$, which contradicts the minimality of $|T|$. Hence we assume  $\theta\in\overline{\varphi}(x)$. Clearly $\alpha \in \varphi(x)$. Let $P_{v_1}:=P_{v_1}(\alpha,\theta,\varphi)$, $P_{x}:=P_{x}(\alpha,\theta,\varphi)$ and $P_z:=P_z(\alpha,\theta,\varphi)$ be $(\alpha, \theta)$-chains containing vertices $v_1$, $x$ and $z$, respectively. By Claim \ref{CLA:OneComp} we have $P_{v_1}=P_{x}$ which is disjoint with $P_z$. Furthermore Lemma \ref{LEM:EXconn} implies that $E(P_{v_1})\supset E_\theta(\partial(T_m))$, together with the assumption that $\alpha\in\Gamma^f(T_0)$, we get $V(P_z)\cap V(T_m)=\emptyset$. Let $\varphi_1=\varphi /P_z$ be the coloring obtained from $\varphi$ by interchanging colors $\alpha$ and $\theta$ along $P_z$. Observe that $\theta$ is only used on $f$ for $E(T-(T_m\cup \partial(T_m)))$ since $\theta\in\overline{\varphi}(x)$, $f$ is colored by $\alpha$ in $\varphi_1$. Clearly $T_m$ is  $(\varphi,\varphi_1)$ stable. By Lemma \ref{LEM:FixConn}, $T_m$ is an ETT of $(G,e,\varphi_1)$, so is $T$. By Claim \ref{CLA:Non0}, let $\gamma \in \Gamma^f(T_0)-(\varphi_1(E(T))\cup \{\theta\})$.  Say $\gamma \in \overline\varphi(v_2)$ for $v_2 \in V(T_0)$. By Claim \ref{CLA:OneComp} the $(\gamma, \theta)$-chain $P^{'}_{v_2}:=P_{v_2}(\gamma,\theta,\varphi_1)$ is the same with $P^{'}_{x}:=P_{x}(\gamma,\theta,\varphi_1)$, hence it is disjoint with $P^{'}_{z}:=P_{z}(\gamma,\theta,\varphi_1)$. Now we consider $T_{zx}$ obtained from $T$ by switching the order of adding vertices $x$ and $z$. Clearly $T_{zx}$ is an ETT of $(G,e,\varphi_1)$ since $f$ is colored by $\alpha$ in $\varphi_1$. Similarly by Claim \ref{CLA:OneComp} the $(\gamma, \theta)$-chain $P^{'}_{v_2}:=P_{v_2}(\gamma,\theta,\varphi_1)$ is the same with $P^{'}_{z}:=P_{z}(\gamma,\theta,\varphi_1)$. Now we reach a contradiction.
\qed

We now prove the following claim which gives a contradiction to {\bf Assumption 1} and  completes the proof of this subcase.
\begin{CLA}\label{CLA:min1}
There is a coloring $\phiv_1\in \mathcal{C}^k(G-e)$ such that $T$ is a non-elementary ETT of $(G, e, \phiv_1)$, $T_m$ is $(\phiv, \phiv_1)$-stable, and
 $|\phibar_1(V(T_0))\cap \phiv_1(E(T)-E(T_0))| > |\phibar(V(T_0))\cap \phiv(E(T)-E(T_0))|$.
 \end{CLA}

 \proof Following the linear order $\prec_{\ell}$, let $e_1$ be the  first edge in $E(T-T_m)$ with $\phiv(e_1) = \alpha$, and let $y_1$ be the in-end of $e_1$.  Pick a missing color $\beta_1\in \phibar(y_1)$.  Note that, since $\phiv(e_1) = \alpha$ and $\alpha \in \Gamma^f(T_0) - \phiv(E(T_m))$, $e_1\notin \partial(T_m)$. Hence  $y_1\in V(T) - V(T_m)$.    Let $P_{v_1}:=P_{v_1}(\alpha,\beta_1,\phiv)$, $P_{y_1}:=P_{y_1}(\alpha,\beta_1,\phiv)$, and $P_z:=P_z(\alpha,\beta_1,\phiv)$ be $(\alpha, \beta_1)$-chains containing $v_1$, $y_1$ and $z$, respectively.  By Claim~\ref{CLA:OneComp},  $P_{v_1} = P_{y_1}$, which in turn shows that it is disjoint from $P_z$. By Lemma~\ref{LEM:EXconn}, $E(P_{v_1})\supset E_{\beta_1}(\partial(T_m))$, so $V(P_z)\cap V(T_m) = \emptyset$.

 Consider the coloring $\phiv_1 = \phiv/P_z$.  Since $V(P_z)\cap V(T_m) = \emptyset$, $T_m$ is $(\phiv, \phiv_1)$-stable. By Lemma~\ref{LEM:FixConn},
 $T_m$ is an ETT of $(G, e, \phiv_1)$.  Since $e_1$ is the first edge colored with $\alpha$ along $\prec_{\ell}$, we have that  $e_1 \prec_{\ell} e^*$ for all edges  $e^*$ colored with $\beta_1$. So, $T$ is an ETT of $(G, e, \phiv_1)$.  Note that $e_1\in E(P_{y_1}) = E(P_{v_1})$, which in turn gives $\phiv_1(e_1) = \alpha$.  We also note that $\beta_1\in \phibar_1(z)\cap \phibar_1(y_1)$.

By Claim~\ref{CLA:Non0}, there is a color $\gamma \in \Gamma^f(T_0)-\phiv(E(T))$. Let $u\in V(T_0)$ such that $\gamma\in \phibar(u)$.  Let $Q_u:= P_u(\gamma, \beta_1, \phiv_1)$, $Q_{y_1}:=P_{y_1}(\gamma, \beta_1, \phiv_1)$ and $Q_z:=P_z (\gamma, \beta_1, \phiv_1)$ be $(\gamma, \beta_1)$-chains containing $u$, $y_1$ and $z$, respectively. By Claim~\ref{CLA:OneComp}, $Q_u = Q_{y_1}$, so $Q_u$ and $Q_z$ are disjoint. By Lemma~\ref{LEM:EXconn}, $E(Q_u)\supset E_{\beta_1}(\partial(T_m))$, so
$V(Q_z)\cap V(T_m) = \emptyset$.  Let $\phiv_2 = \phiv_1/Q_z$ be a coloring obtained from $\phiv_1$ by interchanging colors on $Q_z$.  Since $V(Q_u)\cap V(T_m) = \emptyset$, $T_m$ is an ETT of $(G, e, \phiv_2)$. Since $\gamma\in \phibar(V(T_0)) -\phiv(E(T))$, $T_m$ can be extended to $T$ as an ETT in $\phiv_2$.  Since $\gamma\in \phibar_2(z)\cap \phibar_2(u)$, by Claim~\ref{CLA:a-used}, we have $\gamma\in \phiv_2(E(T-T_m))$. Since $e_1\in Q_{y_1} = Q_u$, the color $\alpha$ assigned to $e_1$ is unchanged.  Thus,
\[
\phibar_2(V(T_0))\cap \phiv_2(E(T)-E(T_0)) \supseteq  ( \phibar(V(T_0))\cap \phiv(E(T)-E(T_0)))\cup\{\gamma\},
\]
 and $\alpha \in \phibar(V(T_0))\cap \phiv(E(T))$.  So, Claim~\ref{CLA:min1} holds. \qed

{\bf \noindent Subcase 2.2:} $y\in V(T_{m})$, i.e. $f\in \partial(T_m)$.

The following two claims are similar to Claims~\ref{CLA:a-used} and \ref{CLA:min1} in Subcase 2.1, which lead to a contradiction to {\bf Assumption 1}.  Their proofs respectively are similar to those of the previous two claims. However, for the completeness, we still give the details.
\begin{CLA}\label{CLA:a-used2}
Color $\alpha$ is used in $E(T-T_m)$, i.e., $\alpha\in \phiv(E(T-T_m))$.
\end{CLA}
\proof Suppose on the contrary $\alpha \notin \phiv(E(T-T_m))$.  By Claim~\ref{CLA:a-free}, we assume that $\alpha \notin \phiv(E(T_m))$,  so  $\alpha \notin \phiv(E(T))$. Let $\phiv(f)=\theta$. As $f\in \partial(T_m)$ is not a connecting edge  and $T_m$ is closed, we know that
there exists $w\in V(T-T_m)$
such that $\theta\in \phibar(w)$.
Consider the $(\alpha,\theta)$-chain $P_{v_1}:=P_{v_1}(\alpha,\theta,\phiv)$.
By Lemma~\ref{LEM:EXconn}, $E(P_{v_1}) \supset E_{\theta}(\partial(T_m))$. So, $f\in E(P_{v_1})$ and $z$ is the other end of $P_{v_1}$.
Then,  $P_{w}:=P_{w}(\alpha,\theta,\phiv)$ is disjoint from $P_{v_1}$, which in turn shows $V(P_{w})\cap V(T_m) = \emptyset$.  Let
$\phiv_1 = \phiv/P_{w}$.    Since $V(P_{w})\cap V(T_m) = \emptyset$,  $T_m$ is $(\phiv, \phiv_1)$-stable. By Lemma~\ref{LEM:FixConn}, $T_m$ is an ETT of $(G, e, \phiv_1)$. Since $\alpha$ is not used in $T-T_m$, $T_m$ can be extended to $T'$ as an ETT of $(G, e, \phiv_1)$. Note that $\alpha\in \phibar_1(v_1)\cap \phibar_1(w)$. So, $T'$ is not elementary, which gives a contradiction to the minimality of $|T|$.   \qed

\begin{CLA}\label{CLA:min2}
There is a coloring $\phiv_1\in \mathcal{C}^k(G-e)$ such that $T$ is a non-elementary ETT of $(G, e, \phiv_1)$, $T_m$ is $(\phiv, \phiv_1)$-stable, and
 $|\phibar_1(V(T_0))\cap \phiv_1(E(T)-E(T_0))| > |\phibar(V(T_0))\cap \phiv(E(T)-E(T_0))|$.
 \end{CLA}

 \proof Following the linear order $\prec_{\ell}$, let $e_1$ be the  first edge in $E(T-T_m)$ with $\phiv(e_1) = \alpha$, and let $y_1$ be the in-end of $e_1$.  Pick a missing color $\beta_1\in \phibar(y_1)$.  Since $\phiv(e_1) = \alpha\in \phibar(V(T_0))$ and $T_m$ is closed, $e_1\notin \partial(T_m)$. Hence,   $y_1\in V(T) - V(T_m)$.    Let $P_{v_1}:=P_{v_1}(\alpha,\beta_1,\phiv)$, $P_{y_1}:=P_{y_1}(\alpha,\beta_1,\phiv)$, and $P_z:=P_z(\alpha,\beta_1,\phiv)$ be $(\alpha, \beta_1)$-chains containing $v_1$, $y_1$ and $z$, respectively.  By Claim~\ref{CLA:OneComp},  $P_{v_1} = P_{y_1}$, which in turn shows that it is disjoint from $P_z$. By Lemma~\ref{LEM:EXconn}, $E(P_{v_1})\supset E_{\beta_1}(\partial(T_m))$, so $V(P_z)\cap V(T_m) = \emptyset$.

 Consider the coloring $\phiv_1 = \phiv/P_z$.  Since $V(P_z)\cap V(T_m) = \emptyset$, $T_m$ is $(\phiv, \phiv_1)$-stable. By Lemma~\ref{LEM:FixConn},
 $T_m$ is an ETT of $(G, e, \phiv_1)$.  Since $e_1$ is the first edge colored with $\alpha$ along $\prec_{\ell}$, we have that  $e_1 \prec_{\ell} e^*$ for all edges  $e^*$ with $\phiv_1(e^*) = \beta_1$.  So, $T$ is an ETT of $(G, e, \phiv_1)$.  Note that $e_1\in E(P_{y_1}) = E(P_{v_1})$, which in turn gives $\phiv_1(e_1) = \alpha$.  We also note that $\beta_1\in\phibar_1(z)\cap \phibar_1(y_1)$.

By Claim~\ref{CLA:Non0}, there is a color $\gamma \in \Gamma^f(T_0)-\phiv(E(T))$. Let $u\in V(T_0)$ such that $\gamma\in \phibar(u)$.  Let $Q_u:= P_u(\gamma, \beta_1, \phiv_1)$, $Q_{y_1}:=P_{y_1}(\gamma, \beta_1, \phiv_1)$ and $Q_z:=P_z (\gamma, \beta_1, \phiv_1)$ be $(\gamma, \beta_1)$-chains containing $u$, $y_1$ and $z$, respectively. By Claim~\ref{CLA:OneComp}, $Q_u = Q_{y_1}$, so $Q_u$ and $Q_z$ are disjoint. By Lemma~\ref{LEM:EXconn}, $E(Q_u)\supset E_{\beta_1}(\partial(T_m))$, so
$V(Q_z)\cap V(T_m) = \emptyset$.  Let $\phiv_2 = \phiv_1/Q_z$ be the  coloring obtained from $\phiv_1$ by interchanging colors on $Q_z$.  Since $V(Q_u)\cap V(T_m) = \emptyset$, $T_m$ is an ETT of $(G, e, \phiv_2)$. Since $\gamma\in \phibar(V(T_0)) -\phiv(E(T))$, $T_m$ can be extended to $T$ as an ETT in $\phiv_2$.  Since $\gamma\in \phibar_2(z)\cap \phibar_2(u)$, by Claim~\ref{CLA:a-used}, we have $\gamma\in \phiv_2(E(T-T_m))$. Since $e_1\in Q_{y_1} = Q_u$, $\varphi_1(e_1) = \varphi(e_1) = \alpha$.  Thus,
\[
\phibar_2(V(T_0))\cap \phiv_2(E(T)-E(T_0)) \supseteq  (\phibar(V(T_0))\cap \phiv(E(T) -E(T_0)))\cup\{\gamma\},
\]
and $\alpha\in \phibar(V(T_0))\cap\phiv(E(T)$.   So, Claim~\ref{CLA:min2} holds. \qed

We now complete the proof of Theorem~\ref{LEM:Elem}. \qed

 Combining Theorem~\ref{LEM:Elem}  and Lemma~\ref{LEM:ETT>F}, we obtain the following result.
 \begin{COR}\label{COR:Main}
 Let $G$ be a $k$-critical graph with $k\ge \Delta+1$.  If $G$ is not elementary, then there is an ETT $T$
 based on  $T_0\in \mathcal{T}$  with $m$-rungs such that $T$ is elementary and
 \[
 |T| \ge |T_0| -2+ \min\{m + |\phibar(V(T_0))|, 2(|\Gamma^f(T_0)|-1) \}.
 \]
 \end{COR}

\section{Proofs of Theorems~\ref{THM:cubic} and
~\ref{THM:Jm19} }
\subsection{Proof of Theorem~\ref{THM:cubic}}
Clearly, we only need to prove Theorem~\ref{THM:cubic} for critical graphs.
\begin{THM}\label{THM:cubic-2}
If $G$ is a $k$-critical graph with $k\geq\Delta+\sqrt[3]{\Delta/2}$,
then $G$ is elementary.
\end{THM}
\proof Suppose on the contrary that $G$ is not elementary. By Corollary~\ref{COR:Main}, let $T$ be an ETT of a $k$-triple $(G, e, \phiv)$ based on $T_0\in \mathcal{T}$
with $m$-rungs  such that $V(T)$ is elementary and
\[
|T| \ge |T_0| -2+ \min\{m + |\phibar(V(T_0))|, 2(|\Gamma^f(T_0)|-1) \}.
\]
Since $m\ge 1$ and $|\phibar(V(T_0))| \ge (k-\Delta) |T_0| +2$, we have
$|T_0| -2 + m +|\phibar(V(T_0))| \ge (k-\Delta +1)|T_0|  +1$.  Following Scheide~\cite{Scheide-2010}, we may assume that $T_0$ is a balanced Tashkinov tree with height $h(T_0) \ge 5$. So, $|\phiv(E(T_0))| \le \frac{|T_0| -1} 2$, which in turn gives
\[
|\Gamma^f(T_0)| = |\phibar(V(T_0))| - |\phiv(E(T_0))| \ge (k-\Delta -\frac 12)|T_0| +\frac{5}{2}.
\]
Hence
\[
|T_0| -2 + 2(|\Gamma^f(T_0)|-1) \ge 2(k -\Delta)|T_0| +1 \ge (k-\Delta +1)|T_0| +1.
\]
Therefore, in any case, we have the following inequality
\begin{equation}\label{EQ:T>}
|T| \ge (k -\Delta +1) |T_0| +1.
\end{equation}
By Corollary~\ref{COR:T-order-exit-vertex}, $|T_0| \ge 2(k-\Delta) +1$. Following ~(\ref{EQ:T>}), we get the inequality below.
\begin{equation}\label{EQ:T>2}
|T| \ge (k-\Delta +1)(2(k-\Delta) +1) + 1= 2(k-\Delta)^2 + 3(k-\Delta) + 2.
\end{equation}
Since $T$ is elementary, we have $k   \ge  |\phibar(V(T))| \ge (k-\Delta)|T|  +2$.  Plugging into (\ref{EQ:T>2}),
we get the following inequality.
\[
k  \ge  2(k-\Delta)^3 + 3(k - \Delta)^2 + 2(k-\Delta) + 2.
\]
Solving the above inequality, we obtain that $k < \Delta +\sqrt[3]{\Delta/2}$, giving a contradiction to
$k\geq\Delta +\lceil \sqrt[3]{\Delta/2}\rceil$. \qed

\subsection{Proofs of Theorem~\ref{THM:Jm19} and Corollary \ref{COR:Delta23}}
We will need the following  observation from~\cite{StiebSTF-Book}. For completeness, we give its proof here.

\begin{LEM}\label{LEM:elmentary-Jm}
Let $s\ge 2$ be a positive integer and $G$ be a $k$-critical graph with $k>   \frac{s}{s-1}\Delta+\frac{s-3}{s-1}$.
For any edge $e\in E(G)$,  if $X\subseteq V(G)$ is an elementary set with respect to
a coloring $\phiv\in \mathcal{C}^k(G-e)$
such that $V(e)\subseteq X$,
then $|X|\le s-1$.
\end{LEM}
\proof Otherwise, assume $|X| \ge s$.   Since $X$ is elementary,
$k \ge |\phibar(X)| \ge (k-\Delta)|X| +2 \ge s(k-\Delta) +2$, which in turn gives
\[
\Delta \ge (s-1)(k-\Delta) +2 > (\Delta +(s-3)) +2 = \Delta + s -1 > \Delta,
\]
a contradiction. \qed

Clearly, to prove Theorem \ref{THM:Jm19}, it is sufficient to restrict our consideration to critical graphs.
\begin{THM}\label{THM:Jm19-2}
If $G$ is  a $k$-critical graph with $k >\frac{23}{22}\Delta+\frac{20}{22}$,
then $G$ is elementary.
\end{THM}
\proof Suppose, on the contrary, $G$ is not elementary.  By Corollary~\ref{COR:Main}, let $T$ be an ETT of a $k$-triple $(G, e, \phiv)$ based on $T_0\in \mathcal{T}$
with $m$-rungs  such that $V(T)$ is elementary and
\[
|T| \ge |T_0| -2+ \min\{m + |\phibar(V(T_0))|, 2(|\Gamma^f(T_0)|-1) \}.
\]

By Lemma~\ref{LEM:elmentary-Jm}, $|T| \le 22$.  We will show that $|T| \ge 23$ to lead a contradiction.
By Lemma~\ref{LEM:Small-T}, we have $|T_0| \ge 11$.
Since $G$ is not elementary, $V(T_0)$ is not strongly closed, so $T\supsetneq T_0$. In particular, we have $m\ge 1$. Since $e\in E(T_0)$, we have $|\phibar(V(T_0))| \ge |T_0| +2$.  Thus,
\begin{equation}\label{EQ:>23}
|T_0| -2+ m + |\phibar(V(T_0))| \ge 2|T_0| +1 \ge 2\times 11 +1 = 23.
\end{equation}

Following Scheide~\cite{Scheide-2010}, we may assume that $T_0$ is a balanced Tashkinov tree with height $h(T_0) \ge 5$, which in turn gives $|\phiv(E(T_0))| \le (|T_0| -1)/2$.  So, $|\Gamma^f(T_0)| \ge |T_0| +2 - (|T_0|-1)/2 \ge (|T_0| +5)/2$.  Thus,
\begin{equation} \label{EQ:>25}
|T_0| -2  +2(|\Gamma^f(T_0)|-1) \ge 2|T_0| +1 \ge 23.
\end{equation}

 Combining (\ref{EQ:>23}) and (\ref{EQ:>25}), we get $|T| \ge 23$, giving a contradiction. \qed

We now give a proof of Corollary \ref{COR:Delta23} and recall that Corollary \ref{COR:Delta23} is stated as follows.

\begin{COR}\label{COR:Delta23-2} If $G$ is a graph with $\Delta\leq 23$ or $|G|\leq 23$, then $\chi'\leq\max\{\Delta+1,\lceil\chi'_f\rceil\}$.
\end{COR}

\proof We assume that $G$ is critical. Otherwise, we prove the corollary for a critical subgraph of $G$ instead. If $\Delta\leq 23$, then $\left\lfloor\frac{23}{22}\Delta+\frac{20}{22}\right\rfloor=\left\lfloor\Delta+\frac{\Delta+20}{22}\right\rfloor\leq \Delta+1$. If $\chi'\leq \Delta+1$, we are done. Otherwise, we assume that $\chi'\geq \Delta+2\geq \frac{23}{22}\Delta+\frac{20}{22}$. By Theorem \ref{THM:Jm19}, we have $\chi'=\lceil \chi_{f}'\rceil$.

Assume that $|G|\leq23$. If $\chi'\leq\Delta+1$, then we are done. Otherwise, $\chi'=k+1$ for some integer $k\geq \Delta+1$. By Corollary \ref{COR:Main}, let $T$ be an ETT of a $k$-triple $(G, e, \phiv)$
based on $T_0\in \mathcal{T}$
with $m$-rungs such that $V(T)$ is elementary and
\[
|T| \ge |T_0| -2+ \min\{m + |\phibar(V(T_0))|, 2(|\Gamma^f(T_0)|-1) \}.
\]

By Lemma~\ref{LEM:Small-T}, we have $|T_0| \ge 11$.
Suppose that $G$ is not elementary, then $V(T_0)$ is not strongly closed, so $T\supsetneq T_0$. In particular, we have $m\ge 1$. Since $e\in E(T_0)$, we have $|\phibar(V(T_0))| \ge |T_0| +2$.  Thus,
\begin{equation}\label{EQ:>23-2}
|T_0| -2+ m + |\phibar(V(T_0))| \ge 2|T_0| +1 \ge 2\times 11 +1 = 23.
\end{equation}

Following Scheide~\cite{Scheide-2010}, we may assume that $T_0$ is a balanced Tashkinov tree with height $h(T_0) \ge 5$, which in turn gives $|\phiv(E(T_0))| \le (|T_0| -1)/2$.  So, $|\Gamma^f(T_0)| \ge |T_0| +2 - (|T_0|-1)/2 \ge (|T_0| +5)/2$.  Thus,
\begin{equation} \label{EQ:>25-2}
|T_0| -2  +2(|\Gamma^f(T_0)|-1) \ge 2|T_0| +1 \ge 23.
\end{equation}

 Combining (\ref{EQ:>23-2}) and (\ref{EQ:>25-2}), we get $|T| \ge 23$. Then $|G|\geq |T|\geq 23$. Therefore, $|G|=23$ and $G$ is elementary, giving a contradiction. \qed

 \section{Acknowledgement}

  We thank Guangming Jing for comments that greatly improved the manuscript.

  \bibliographystyle{plain}
\bibliography{Gold2016Dec11}

\begin{thebibliography}{10}

\bibitem{Andersen}
Lars~D{\o}vling Andersen.
\newblock On edge-colourings of graphs.
\newblock {\em Math. Scand.}, 40(2):161--175, 1977.

\bibitem{CYZ-2011}
Guantao Chen, Xingxing Yu, and Wen'an Zang.
\newblock Approximating the chromatic index of multigraphs.
\newblock {\em J. Comb. Optim.}, 21(2):219--246, 2011.

\bibitem{Edmonds65}
Jack Edmonds.
\newblock Maximum matching and a polyhedron with {$0,1$}-vertices.
\newblock {\em J. Res. Nat. Bur. Standards Sect. B}, 69B:125--130, 1965.

\bibitem{FavrST06}
Lene~Monrad Favrholdt, Michael Stiebitz, and Bjarne Toft.
\newblock Graph edge coloring: Vizing's theorem and goldberg's conjecture.
\newblock Preprint DMF-2006-10-003, IMADA-PP-2006-20, University of Southern
  Demark, 2006.

\bibitem{Goldberg}
Mark~K. Goldberg.
\newblock On multigraphs of almost maximal chromatic class\,(russian).
\newblock {\em Discret. Analiz.}, 23:3--7, 1973.

\bibitem{Goldberg-1984}
Mark~K. Goldberg.
\newblock Edge-coloring of multigraphs: recoloring technique.
\newblock {\em J. Graph Theory}, 8(1):123--137, 1984.

\bibitem{Gupta67}
Ram~Prakash Gupta.
\newblock {\em Studies in the Theory of Graphs}.
\newblock 1967.
\newblock Thesis (Ph.D.)--Tata Institute of Fundamental Research, Bombay.

\bibitem{Holyer81}
Ian Holyer.
\newblock The {NP}-completeness of edge-coloring.
\newblock {\em SIAM J. Comput.}, 10(4):718--720, 1981.

\bibitem{Jakobsen73}
Ivan~Tafteberg Jakobsen.
\newblock Some remarks on the chromatic index of a graph.
\newblock {\em Arch. Math. (Basel)}, 24:440--448, 1973.

\bibitem{Kahn96}
Jeff Kahn.
\newblock Asymptotics of the chromatic index for multigraphs.
\newblock {\em J. Combin. Theory Ser. B}, 68(2):233--254, 1996.

\bibitem{Kierstead84}
Henry~Andrew Kierstead.
\newblock On the chromatic index of multigraphs without large triangles.
\newblock {\em J. Combin. Theory Ser. B}, 36(2):156--160, 1984.

\bibitem{McDonaldSurvey15}
Jessica McDonald.
\newblock Edge-colourings.
\newblock In L.~W. Beineke and R.J. Wilson, editors, {\em Topics in topological
  graph theory}, pages 94--113. Cambridge University Press, 2015.

\bibitem{Nishizeki-Kashiwagi-1990}
Takao Nishizeki and Kenichi Kashiwagi.
\newblock On the {$1.1$} edge-coloring of multigraphs.
\newblock {\em SIAM J. Discrete Math.}, 3(3):391--410, 1990.

\bibitem{Scheide-2010}
Diego Scheide.
\newblock Graph edge colouring: {T}ashkinov trees and {G}oldberg's conjecture.
\newblock {\em J. Combin. Theory Ser. B}, 100(1):68--96, 2010.

\bibitem{Seymour}
Paul Seymour.
\newblock On multicolourings of cubic graphs, and conjectures of {F}ulkerson
  and {T}utte.
\newblock {\em Proc. London Math. Soc. (3)}, 38(3):423--460, 1979.

\bibitem{Shannon49}
Claude~Elwood Shannon.
\newblock A theorem on coloring the lines of a network.
\newblock {\em J. Math. Physics}, 28:148--151, 1949.

\bibitem{StiebSTF-Book}
Michael Stiebitz, Diego Scheide, Bjarne Toft, and Lene~Monrad Favrholdt.
\newblock {\em Graph edge coloring}.
\newblock Wiley Series in Discrete Mathematics and Optimization. John Wiley \&
  Sons, Inc., Hoboken, NJ, 2012.
\newblock Vizing's theorem and Goldberg's conjecture, With a preface by
  Stiebitz and Toft.

\bibitem{Tashkinov-2000}
Vladimir~Aleksandrovich Tashkinov.
\newblock On an algorithm for the edge coloring of multigraphs.
\newblock {\em Diskretn. Anal. Issled. Oper. Ser. 1}, 7(3):72--85, 100, 2000.

\bibitem{Vizing64}
Vadim~Georgievich Vizing.
\newblock On an estimate of the chromatic class of a {$p$}-graph.
\newblock {\em Diskret. Analiz No.}, 3:25--30, 1964.

\end{thebibliography}





\bibliographystyle{amsplain}

\end{document}